\documentclass[12pt]{amsart}

\usepackage[margin=1.15in]{geometry}

\usepackage{amsmath,amscd,amssymb,amsfonts,latexsym,wasysym, mathrsfs, mathtools,hhline,color}
\usepackage[all, cmtip]{xy}
\usepackage{csquotes}
\usepackage{url}

\definecolor{hot}{RGB}{65,105,225}

\usepackage[pagebackref=true,colorlinks=true, linkcolor=hot ,  citecolor=hot, urlcolor=hot]{hyperref}
\usepackage{ textcomp }
\usepackage{ tipa }

\usepackage{graphicx,enumerate}

\theoremstyle{plain}
\newtheorem{theorem}{Theorem}[]
\newtheorem{prop}[theorem]{Proposition}

\newtheorem{lm}[theorem]{Lemma}

\newtheorem{cor}[theorem]{Corollary}
\newtheorem{conj}[theorem]{Conjecture}
\newtheorem{lemma}[theorem]{Lemma}
\newtheorem{thrm}[theorem]{Theorem}

\theoremstyle{definition}

\newtheorem{defn}[theorem]{Definition}
\newtheorem{question}[theorem]{Question}
\newtheorem{rmk}[theorem]{Remark}

\newtheorem{ex}[theorem]{Example}
\newtheorem*{ex*}{Example}

\def\be{\begin{equation}}
\def\ee{\end{equation}}

\def\bt{\begin{thrm}}
\def\et{\end{thrm}}

\def\bc{\begin{cor}}
\def\ec{\end{cor}}

\def\br{\begin{rmk}}
\def\er{\end{rmk}}

\def\bp{\begin{prop}}
\def\ep{\end{prop}}

\def\bl{\begin{lm}}
\def\el{\end{lm}}

\def\bex{\begin{ex}}
\def\eex{\end{ex}}

\def\bd{\begin{defn}}
\def\ed{\end{defn}}

\newcommand\sP{{\mathcal P}}

\newcommand\sF{{\mathcal F}}

\def\bK{\mathbf{K}}

\newcommand\K{{\mathbb{K}}}
\newcommand\bP{{\mathbb{P}}}

\newcommand\Z{{\mathbb{Z}}}

\newcommand\C{{\mathbb{C}}}

\DeclareMathOperator{\codim}{codim}              
                  
\DeclareMathOperator{\id}{id}                    

\def\bC{\mathbb{C}}

\def\bP{\mathbb{P}}

\def\bZ{\mathbb{Z}}

\def\bK{\mathbb{K}}

\title[A new generic vanishing theorem on homogeneous varieties]{A new generic vanishing theorem on homogeneous varieties and the positivity conjecture  for triple intersections of Schubert cells}

\author[J. Sch\"urmann ]{J\"org Sch\"urmann}
\address{J.  Sch\"urmann : Mathematische Institut,
          Universit\"at M\"unster,
          Einsteinstr. 62, 48149 M\"unster,
          Germany.}
\email {jschuerm@uni-muenster.de}

\author[C. Simpson ]{Connor Simpson}
\address{C. Simpson : Department of Mathematics, University of Wisconsin-Madison, 480 Lincoln Drive, Madison WI 53706-1388, USA}
\email {csimpson6@wisc.edu}

\author[B. Wang ]{Botong Wang}
\address{B. Wang : Department of Mathematics, University of Wisconsin-Madison, 480 Lincoln Drive, Madison WI 53706-1388, USA}
\email {wang@math.wisc.edu}


\thanks{J. Sch\"urmann is funded by the Deutsche Forschungsgemeinschaft (DFG, German Research Foundation) Project-ID 427320536 -- SFB 1442, as well as under Germany's Excellence Strategy EXC 2044 390685587, Mathematics M\"unster: Dynamics--Geometry--Structure. C. Simpson thanks Xuhua He for his hospitality and support during Spring 2023. B. Wang is partially supported by a Sloan fellowship. The authors also thank the University of Wisconsin-Madison for funding our collaboration.}

\begin{document}

\maketitle
\begin{abstract} In this paper we prove a new generic vanishing theorem for $X$  a complete homogeneous variety with respect to an action of a connected algebraic group. Let $A, B_0\subset X$ be locally closed affine subvarieties, and assume that $B_0$ is smooth and pure dimensional. Let $\sP$ be a perverse sheaf on $A$ and let $B=g B_0$ be a generic translate of $B_0$. Then our theorem implies
$(-1)^{\codim B}\chi(A\cap B, \sP|_{A\cap B})\geq 0$. As an application, we prove in full generality
a positivity conjecture about the signed Euler characteristic of generic  triple intersections of Schubert cells. Such Euler characteristics are known to be the structure constants for the multiplication of the Segre-Schwartz-MacPherson classes of these Schubert cells.
\end{abstract}

\section{Introduction}
Let $G$ be a complex semisimple, simply connected, linear algebraic group, and fix a Borel subgroup $B$ with a maximal torus $T\subseteq B$. Let $B^-$ denote the opposite Borel subgroup, with $P\subset G$ a parabolic subgroup containing $B$. Let $W$ be the Weyl group of $B$, 
and let $w_0\in W$ be the longest element in $W$; then $B^- = w_0 B w_0$. 
Let $W_P$ be the subgroup of $W$ generated by the simple reflections
in $P$ and denote by $W^P\subset W$ the set of minimal length representatives
for the cosets of $W_P$ in $W$.

Given any element $\lambda\in W^P$, let $X_\lambda^\circ$ and $X^{\lambda\circ}$ be the corresponding Schubert cell and the opposite Schubert cell in the partial flag variety $X:=G/P$ respectively. Denote by $X_\lambda$ and $X^\lambda$ the corresponding Schubert variety and opposite Schubert variety. The Schubert cells $X_\lambda^\circ$ (resp. opposite Schubert cell $X^{\lambda\circ}$) are $B$-orbits (resp. $B^-$-orbits) for the left multiplication action on $X$. Moreover, 
\[
X_\lambda=\bigsqcup_{\xi\in W^P,\, \xi\leq \lambda}X_{\xi}^\circ\quad \text{and}\quad X^\nu=\bigsqcup_{\xi\in W^P,\, \xi\geq \nu}X^{\xi\circ}
\]
are algebraic Whitney stratifications of the Schubert variety $X_\lambda$ and the opposite Schubert variety $X^\nu$, respectively. 
By a classical result of Richardson \cite{Ric92}, these stratifications 
of $X_\lambda$ and $X^{\nu}$ are transversal in $X$, so that the Richardson variety $X_\lambda \cap X^{\nu}$ admits an induced algebraic Whitney stratification.


For a locally closed algebraic subset $U\subset X$, let 
$$c_{SM}(U):=c_M(1_U)\quad \text{and} \quad s_{SM}(U):=\frac{c_{SM}(U)}{c(TX)}\in H^*(X,\bZ)$$
be the  Chern-Schwartz-MacPherson (CSM) class and the Segre-Schwartz-MacPherson (SSM) class of $U$ or $1_U$ in the ambient variety $X$, respectively. Then by an intersection formula of the first author \cite[Theorem 1.2]{Sch17}, the transversality of these stratifications implies the formula
\begin{equation*}
c_{SM}(X_\lambda^\circ)\cdot s_{SM}(X^{\nu\circ})=c_{SM}(X_\lambda^\circ \cap X^{\nu\circ}) \in H^*(X,\bZ)\:,
\end{equation*}
with  $X_\lambda^\circ \cap X^{\nu\circ}$ an open   Richardson variety.
By counting $T$-fixed points in $X_\lambda^\circ \cap X^{\nu\circ}$ and the functoriality of the MacPherson-Chern class transformation $c_M$, the above formula implies the geometric orthogonality relation (see \cite[Theorem 7.1]{AMSS1}):
\begin{equation}\label{ortho}
\langle c_{SM}(X_\lambda^\circ),  s_{SM}(X^{\nu\circ})\rangle :=\int_X\,c_{SM}(X_\lambda^\circ)\cdot s_{SM}(X^{\nu\circ})
=\chi(X_\lambda^\circ \cap X^{\nu\circ})=\delta_{\lambda, \nu}\:.
\end{equation}
For any $\lambda, \mu, \nu\in W^P$, we introduce the SSM structure constants $a^\nu_{\lambda,\mu}\in \bZ$ via
\begin{equation*}
s_{SM}(X_\lambda^\circ)\cdot s_{SM}(X_\mu^\circ) =\sum_\nu  a^\nu_{\lambda,\mu}\cdot s_{SM}(X_\nu^\circ) \in H^*(X,\bZ)\:.
\end{equation*}

By Kleiman's transversality theorem \cite{Kleiman}, for generic $g\in G$,
 the translate $gX_\mu$ is stratified transversal to the Richardson variety $X_\lambda \cap X^{\nu}$ (with its induced algebraic Whitney stratification).
Then the intersection formula \cite[Theorem 1.2]{Sch17} together with the geometric  orthogonality relation (\ref{ortho}) implies 
(see e.g. \cite{Su21, Kumar, AMSS}):
\begin{equation*}\label{structure}
a^\nu_{\lambda,\mu}=\chi(X_\lambda^\circ\cap gX_\mu^\circ \cap X^{\nu\circ})=\chi(X_\lambda^\circ\cap gX_\mu^\circ \cap hX_{\nu'}^\circ)
\end{equation*}
for $g,h\in G$ generic and $\nu'\cdot W_P=w_0\nu\cdot W_P$. Note that $c_{SM}(gX_\mu^\circ)=c_{SM}(X_\mu^\circ)$ by functoriality of the  MacPherson-Chern class transformation,
together with the fact that the connected group $G$ acts trivially on $H^*(X,\bZ)$. Similarly $c_{SM}(X^{\nu\circ})=c_{SM}(hX_{\nu'}^\circ)$,
 since there is an element $n_0\in G$ (lifting $w_0\in W$) with $n_0X_{\nu'}^\circ=X^{\nu\circ}$ (see e.g. \cite[Eq. (29)]{AMSS1}, also for the corresponding $T$-equivariant result).
Several authors formulated the following conjecture, posed from 2019 onward in talks of L. Mihalcea and A. Knutson, then in \cite{KZJ} and finally in \cite[Conjecture D]{Kumar}, as surveyed in \cite[Conjecture 3]{AMSS}:
\begin{conj}[Signed Euler characteristic of the intersection of three Schubert cells]\label{mainconj}
$$E_{\lambda, \mu, \nu'}:= (-1)^{d}\cdot a^\nu_{\lambda,\mu}= (-1)^{d}\cdot   \chi(X_\lambda^\circ\cap g X_\mu^\circ\cap h X_{\nu'}^\circ)
\geq 0$$
for $g,h\in G$ generic, with $\nu'\cdot W_P=w_0\nu\cdot W_P$ and $d:=\dim (X_\lambda^\circ\cap g X_\mu^\circ\cap h X_{\nu'}^\circ)$. 
\end{conj}

\begin{rmk}
  By \cite[Theorem 2.5]{Su21} or \cite[Eq. 36]{AMSS1}, up to a sign, $E_{\lambda, \mu, \nu'}$ are also the corresponding structure constants for the CSM classes of Schubert cells in a full flag variety $G/B$. 
\end{rmk}

Kumar \cite[Theorem 16]{Kumar} showed for the case of full flag varieties $G/B$, that this conjecture would follow from
 positivity of the Chern-Schwartz-MacPherson classes $c_{SM}(X_\lambda^\circ \cap X^{\nu\circ})$ of all open Richardson varieties in $G/B$, posited by \cite[Conjecture 5]{Kumar} and \cite[Conjecture 9.2]{FGX22}.
In \cite[Theorem 3]{KZJ}, Knutson and Zinn-Justin proved for the case of the $r$-step 
partial flag variety of type $A$, with  $r\leq 4$, that the Euler characteristic
\[
E_{\lambda, \mu, \nu'} := (-1)^{d}\cdot   a^\nu_{\lambda,\mu}= (-1)^{d}\cdot   \chi(X_\lambda^\circ\cap g X_\mu^\circ\cap h X_{\nu'}^\circ)
\]
is equal to a weighted count of certain combinatorial puzzles. In particular, when $r \leq 3$, they obtained the following theorem (see also the recent ICM talk of Knutson 
\cite[Theorem 4]{Kn22}):
\begin{theorem}[\cite{KZJ}]\label{thm_K}
If $X_\lambda^\circ$, $X_\mu^\circ$, and $X_{\nu'}^\circ$ are three Schubert cells in an $r$-step partial flag variety of type $A$, with $r \leq 3$, then $E_{\lambda, \mu, \nu'}\geq 0$.
\end{theorem}

In this paper, we prove the following new generic vanishing theorem on homogeneous varieties, from which we deduce the Conjecture \ref{mainconj} in full generality for all partial flag varieties $X=G/P$ (with  functors like $f_!,f_*$  already denoting their derived versions).

\begin{theorem}[New generic vanishing theorem]\label{thm_GV}
Let $X$ be a complete homogeneous variety with respect to an action of a connected algebraic group $G'$. Let $A, B_0\subset X$ be locally closed affine subvarieties, and assume that $B_0$ is also smooth and pure dimensional. Let $\sP$ be a perverse sheaf on $A$ and $B=g B_0$ be a generic translate of $B_0$.
Then,
\[
H_c^i\big(B, j_{A*}(\sP|_{A\cap B})\big)\cong H^i\big(A, j_{B!}(\sP|_{A\cap B})\big)
=0\quad \text{for all}\quad i\neq -\codim_X B
\]
where $j_B: A\cap B\to A$ and $j_A: A\cap B\to B$ the inclusion maps.
\end{theorem}

Note that such a complete homogeneous variety $X$ is isomorphic to a product $X\cong Ab\times G/P$ of an abelian variety $Ab$ and a partial flag variety 
$G/P$ as before (see, e.g. \cite[Theorem 2.6]{Bri12}). 

A straightforward consequence of  Theorem \ref{thm_GV} is the following. 
\begin{cor}\label{cor_Euler}
Under the notations of the theorem, 
\[
(-1)^{\codim B}\chi(A\cap B, \sP|_{A\cap B})\geq 0. 
\]
\end{cor}
Here we assume that $\sP$ is a perverse sheaf of vector spaces, with finite dimensional stalks (over a given field $\bK$, e.g. $\bK=\bC$, in which case $\sP$ is isomorphic the de Rham  complex of a holonomic D-module), so that the Euler characteristics
$$\chi\left(H_c^i\big(B, j_{A*}(\sP|_{A\cap B})\big)\right) \quad \text{and} \quad \chi\left(H^i\big(A, j_{B!}(\sP|_{A\cap B})\big)\right)$$
only depend on the associated constructible function $\chi_{stalk}$ (given by the stalkwise Euler characteristic):
$$\chi_{stalk}\left(j_{A*}(\sP|_{A\cap B})\right) = j_{A*}\chi_{stalk}(\sP|_{A\cap B}) 
 \quad \text{and} \quad 
\chi_{stalk}\left(j_{B!}(\sP|_{A\cap B})\right) = j_{B!}\chi_{stalk}(\sP|_{A\cap B}) \:,
$$
with $f_!=f_*$ as covariant  functors on the level of constructible functions for an algebraic morphism $f$
(see  \cite[Equation (10.3) and Example 10.4.37]{MS} and \cite[Section 2.3
and Section 6.0.6]{Sch}).
So when $f=j$ is a locally closed inclusion, $j_!=j_*$ on the level of constructible functions is just the extension by zero.\\

\begin{rmk}\label{rem-cs1}
Using composition of functors, we have
\[
H_c^i(B, j_{A!}(\sP|_{A\cap B}))=H_c^i(A\cap B,\sP|_{A\cap B})
\]
and
\[H^i(A, j_{B*}(\sP|_{A\cap B}))
=H^i(A\cap B,\sP|_{A\cap B}). 
\]
We do not claim that any of them is concentrated in one degree.
Our result implies that 
\[
\chi_c(B,j_{A!}(\sP|_{A\cap B}))=\chi(A\cap B, \sP|_{A\cap B})=\chi(A, j_{B*}(\sP|_{A\cap B}))
\]
 has the same sign as $(-1)^{\codim B}$,
by the above discussions.
\end{rmk}

The Conjecture \ref{mainconj} can be deduced from this Corollary \ref{cor_Euler}, with $X=G/P$ and $G'=G$ and the following choices: $A=X_{\lambda}^\circ\cap g X_{\mu}^\circ$, $B_0=X_{\nu'}^\circ$ and $\sP=\bK_{A}[\dim  A]$. 
The results of Theorem \ref{thm_GV} and Corollary \ref{cor_Euler} can easily be inductively applied for $n$ ($n\geq 2$) generic intersections of Schubert cells. 
Let $A=X_{\mu_1}^\circ\cap g_2 X_{\mu_2}^\circ \cap \cdots \cap g_{n-1}X_{\mu_{n-1}}^\circ$, $B_0=X_{\mu_{n}}^\circ$ and 
 $\sP=\bK_{A}[\dim A]$, where $g_2, \ldots, g_{n-1}\in G$ are generic.

\begin{theorem}\label{thm_2}
Let $X_{\mu_1}^\circ, \dots, X_{\mu_{n}}^\circ$ be $n$  Schubert cells ($n\geq 2$) in a flag variety $X=G/P$. Then for  generic $g_2,\dots,g_{n}\in G$ one has 
\[
(-1)^{d}\chi(X_{\mu_1}^\circ\cap g_2 X_{\mu_2}^\circ\cap \cdots \cap g_{n}X_{\mu_{n}}^\circ)\geq 0\:,
\]
with $d=\dim (X_{\mu_1}^\circ\cap g_2 X_{\mu_2}^\circ\cap \cdots \cap g_{n}X_{\mu_{n}}^\circ)$.
\end{theorem}

Note that we need to intersect at least two 
generically translated affine Schubert cells, so that our method applies. 
Of course, $(-1)^{\dim X_{\mu}}\cdot \chi(X_{\mu}^\circ)=(-1)^{\dim X_{\mu}}$ is not always non-negative for one Schubert cell $X_{\mu}^\circ$. 
These iterated open intersections are the open counterparts of the intersection varieties studied in \cite{BC}.
In particular, when $n=2$, one gets $(-1)^{d}\chi(X_{\lambda}^\circ\cap  gX_{\mu}^\circ)\geq 0$, with $d =\dim (X_{\lambda}^\circ \cap 
gX_{\mu}^\circ)$ and $g\in G$ generic. 
This is consistent with the well-known result that for generic $g$, $X_{\lambda}^\circ \cap 
gX_{\mu}^\circ$ is isomorphic to an open Richardson variety $X_\lambda^\circ \cap X^{\nu\circ}$  (see \cite[Remark 2.2]{BC}), so that this signed Euler characteristic is one for $d=0$ and zero otherwise (by the geometric orthogonality relation (\ref{ortho})).\\

In this paper, Theorem \ref{thm_GV} is deduced from Artin's vanishing theorem for perverse sheaves (see  \cite[Theorem 10.3.59]{MS} and 
\cite[Corollary 6.0.4]{Sch}),
Kleiman's transversality theorem \cite{Kleiman}, 
 together with the following key technical proposition:
\begin{prop}[Generic base change isomorphism]\label{prop_diag}
Let $X$ be a  homogeneous variety with respect to an action of a connected algebraic group $G'$. Let $A, B_0$ be locally closed subvarieties of $X$. Let $\sF$ be a constructible complex on $A$ . Let $B=g B_0$ be a generic translate of $B_0$. Denote various inclusion maps as in the diagram
\begin{equation*}
\xymatrix{
A\cap B\ar[d]^{j=j_A}\ar[rr]^{\quad j_B}&& A\ar[d]^{i_A}\\
B\ar[rr]^{i_B}&&X.
}
\end{equation*}
Then there exists a canonical base change isomorphism
in $D_c^b(X)$,
\begin{equation}\label{basech1}
i_{B!}j_{A*}(\sF|_{A\cap B})\cong i_{A*}j_{B!}(\sF|_{A\cap B}) \:.
\end{equation}
\end{prop}

\noindent\textbf{Conventions.} Throughout this paper, we fix a commutative $\bK$ ring with unit as the coefficient ring, so that every (bounded) constructible complex and perverse sheaf has coefficients over $\bK$. When taking Euler characteristics, we assume $\bK$ is a field and all stalks of our constructible sheaf complexes are finite dimensional over $\bK$.
All the pushforward and pullback functors are derived functors. We work over complex algebraic varieties, and all stratifications are assumed to be algebraic
Whitney stratifications. Given finitely many bounded constructible complexes on a given complex variety, there is always such a stratfication with all the sheaf complexes constructible with respect to this stratification, i.e. all their cohomology sheaves are locally constant on all strata.

\section{Motivation of Theorem \ref{thm_GV}}
In the introduction we applied 
our main Theorem \ref{thm_GV} only to the case of a partial flag variety $G/P$, although it applies to a homogenous variety of product type 
$X=Ab\times G/P$,
with $Ab$ an abelian variety. Now, we explain how it
is partially motivated by  generic vanishing theorems for (semi-)abelian varieties. Here a semiabelian variety $N$ fits into an algebraic group extension
$0\to (\bC^*)^n \to N \to Ab\to 0$ of an abelian variety by a complex torus.
The generic vanishing theorem was first developed in the coherent setting by Green and Lazarsfeld \cite{GL}, and it was adapted to perverse sheaves and D-modules in \cite{GaL, Schnell, LMW} and other works. One particular generic vanishing theorem for perverse sheaves is the following. 
\begin{theorem}[\cite{LMW}]
Let $N$ be a semiabelian variety, and let $\sP$ be a perverse sheaf on $N$. Then for a general rank one local system $L$ on $N$, 
\[
H^i(N, \sP\otimes L)=0 \quad\text{for any}\quad i\neq 0. 
\]
\end{theorem}
\begin{cor}\label{cor_SA}
Let $\sP$ be a perverse sheaf on a semiabelian variety $N$. Then $\chi(N, \sP)\geq 0$. 
\end{cor}
\begin{proof}
Since taking tensor product with a rank one local system does not change Euler characteristics, 
\begin{equation*}\label{semiab}
\chi(N, \sP)=\chi(N, \sP\otimes L)\geq 0,
\end{equation*}
where $L$ is a general rank one local system on $N$. 
\end{proof}

The non-negativity of the  Euler characteristic of perverse sheaves on a semiabelian variety was also observed by Franecki and Kapranov using the corresponding  effective characteristic cycles and Kashiwara's index theorem \cite{FK} (see  also \cite[Proposition 8.4]{AMSS0} for a more general version for suitable varieties mapping to an abelian variety $Ab$).
 In the case of an abelian variety it more generally holds for a constructible function $\varphi$ with an effective characteristic cycle, i.e. a non-negative  linear combination of signed Euler obstructions $(-1)^{\dim Z}\cdot Eu_Z$ for $Z\subset Ab$ a closed irreducible subvariety. Moreover the corresponding  Euler characteristic result $\chi(\varphi)\geq 0$ is only a shadow of the fact that in this case the signed MacPherson Chern  class
$c_M^{\vee}(\varphi)$ is effective \cite{AMSS0, ST}. Here $c_M^{\vee}(\varphi)$ differs from $c_M(\varphi)$ by the sign $(-1)^k$ in homological degree $2k$ (e.g. they have the same degree zero part). 

\begin{rmk}
In our new generic vanishing Theorem \ref{thm_GV} we do not use twisting by a generic rank one local system. 
Instead we use in particular the affineness assumptions together with a generic translation on the given ambient complete homogeneous variety, in such a way that 
the underlying constructible functions agree in the sense of Remark~\ref{rem-cs1}.
\end{rmk}

Let us illustrate Corollary  \ref{cor_Euler} in the simplest case of an abelian (or partial flag)  variety given by an elliptic curve  $E$ (or by $\bP^1$).
\begin{ex}
Let $X=E$ be an elliptic curve (or $X=\bP^1$)  with $A=X\backslash S$ and $B=X\backslash S'$ for two finite non-empty subsets $S,S'\subset X$,
with $|S\cup S'|\geq 2$ (as it would be for a generic translate),
so that $A$ and $B$ are smooth and affine of codimension zero. Let the perverse sheaf $\sP=L[1]$ on $A$ be given by a shifted local system $L$ of rank $r\geq 0$. Then 
\[
\chi(A\cap B,\sP)=-r\cdot \big(\chi(X)-|S\cup S'|\big)\geq r\cdot (|S\cup S'|-2) \geq 0,
\]
since $\chi(E)=0$ (and $\chi(\bP^1)$=2).
\end{ex}

The affine space analog of the work of Franecki-Kapranov was studied using stratified Morse theory in \cite{ST, STV}. 
\begin{theorem}[\cite{ST}, \cite{STV}]\label{thm_affine}
Let $\sP$ be a perverse sheaf on the affine space $\C^n$. Let $H$ be a general affine hyperplane in $\C^n$, and let $U=\C^n\setminus H$. Then 
\[
\chi(U, \sP|_{U})\geq 0.
\]
\end{theorem}

Applying our main Theorem \ref{thm_GV} to the case when the ambient complete homogeneous variety is $X=\bP^n$ and $A=B_0=\C^n$, gives a new proof of Theorem \ref{thm_affine}. 

\begin{rmk}\label{rem-FK}
Theorem \ref{thm_GV} implies also directly  the following counterpart for $X$ a complete homogeneous variety with respect to an action of a connected algebraic group $G'$. Let $D_1, D_2\subset X$ be two ample divisors on $X$,
with $A=X\setminus D_1$ and $B_0=X\backslash D_2$ the affine open complements. Let $\sP$ be a perverse sheaf on  $A$. Then
\begin{equation*}
\chi(U, \sP|_{U})\geq 0.
\end{equation*}
for $U:=A\backslash gD_2$ the open complement in $A$ of a translate $gD_2$ of $D_2$ by a generic element $g\in G'$.
\end{rmk}

When the perverse sheaf in Corollary \ref{cor_SA} and Theorem \ref{thm_affine} is of the form $\C_Y[\dim Y]$, where $Y$ is a smooth subvariety of an affine torus or an affine space, there are explicit geometric interpretations of the Euler characteristics (see \cite{Huh, RW}, \cite[Equation (2)]{STV} and \cite[Theorem~1.2]{ST}).\footnote{If $Y$ is a closed smooth subvariety of a semiabelian variety, then $(-1)^{\dim Y}\chi(Y)$ counts the number of degenerate point of $\eta|_{Y}$, where $\eta$ is a general parallel one form on the semiabelian variety (\cite{LMW2}). When $Y$ is singular, then the number of critical points in the smooth locus counts the signed Euler characteristic of the Euler obstruction function $Eu_Y$.} For a smooth subvariety $Y$ of $(\C^*)^n$, $(-1)^{\dim Y}\chi(Y)$ is equal to the number of critical points of $\prod_{1\leq i\leq n}z_i^{u_i}|_Y$, where $z_i$ are the coordinates of $(\C^*)^n$ and $u_i\in \Z$ are general. Similarly, for a smooth subvariety $Y\subset \C^n$, $(-1)^{\dim Y}\chi(Y\setminus H)$ is equal to the number of critical points of $l|_Y$, where $l$ is a general affine linear function and $H=\{l=0\}$. These observations lead us to ask the following question. 
\begin{question}
Let $A$ and $B_0$ be two locally closed and pure dimensional smooth affine subvarieties of a complete  homogeneous variety $X$, 
so that $(-1)^{\dim A\cap B}\chi(A\cap B)\geq 0$
for  $B=g B_0$  a generic translate (by Theorem \ref{thm_GV}). Does $(-1)^{\dim A\cap B}\chi(A\cap B)$ count anything?
\end{question}

The solution of the above question potentially will lead a solution to the following question, which is a constructible function analog of Corollary \ref{cor_Euler}. 
\begin{question}
Let $A, B_0$ be locally closed affine  irreducible subvarieties of a complete homogeneous variety $X$, and assume that $B_0$ is smooth. 
Let $B=g B_0$ be a generic translate. Does the inequality
\[
(-1)^{\dim A\cap B}\chi(Eu_{A\cap B})\geq 0
\]
always hold, where $Eu_{A\cap B}$ denotes the Euler obstruction function of $A\cap B$? 
\end{question}

The result of Theorem \ref{thm_affine} holds for the signed Euler obstruction $(-1)^{\dim Z}\cdot Eu_Z$ of an irreducible algebraic subset 
$Z\subset \bC^n$, since their proof in \cite{ST, STV} is given in terms of stratified Morse theory for constructible functions. However, 
Theorem~\ref{thm_GV}  depends essentially on the base change isomorphism (\ref{basech1}) in the context of constructible sheaf complexes,
which has no counterpart for constructible functions.

\section{Proof of Theorem \ref{thm_GV}}
\subsection{A transversality result}
We first go over some basic transversality results in the subsection to prepare for the proof of Proposition \ref{prop_diag}. 
\begin{lemma}
Let $Z_1$ and $Z_2$ be two locally closed submanifolds of a complex manifold $M$. Then $Z_1$ and $Z_2$ intersect transversally if and only if $Z_1\times Z_2$ intersects the diagonal $\Delta\subset M\times M$ transversally. 
\end{lemma}
\begin{proof}
The two submanifolds $Z_1$ and $Z_2$ intersects transversally at a point $P\in Z_1\cap Z_2$ if and only if the tangent spaces $T_P Z_1$ and $T_P Z_2$ intersect transversally at $P$, which is further equivalent to $T_{(P, P)}Z_1\times Z_2$ intersects transversally with $T_{(P, P)}\Delta$ in $T_{(P, P)}M\times M$ (see also \cite[p.257]{Sch} for a discussion of transversal intersections of product stratifications with respect to suitable regularity conditions like the Whitney condition).
\end{proof}
The following is the well-known transversality theorem of Kleiman. 
\begin{theorem}[\cite{Kleiman}]
Let $X$ be a  homogeneous variety with respect to an action of a connected algebraic group $G'$. Let $\mathfrak{S}=\bigsqcup_{k\in \Lambda_1}S_k$ and $\mathfrak{T}=\bigsqcup_{k'\in \Lambda_2}T_{k'}$ be two Whitney stratifications of $X$ into locally closed smooth complex algebraic subvarieties. For any $g\in G'$, let $g\cdot \mathfrak{T}=\bigsqcup_{k'\in \Lambda_2}g T_{k'}$ be the translate of the stratification $\mathfrak{T}$ by $g$. Then any stratum in $\mathfrak{S}$ and any stratum in $g\cdot \mathfrak{T}$ intersect transversally for $g\in G'$ generic.
\end{theorem}
Combining the above Lemma and Theorem, we have the following generic transversality result. 
\begin{cor}\label{cor_transversal}
Under the notations of the above Lemma, the diagonal $\Delta\subset M\times M$ intersects the stratification $\mathfrak{S}\times g\cdot \mathfrak{T}$ of $M\times M$ transversally  for $g\in G'$ generic.
\end{cor}
\subsection{The proofs}
We first prove Proposition \ref{prop_diag}, then deduce the vanishing Theorem \ref{thm_GV}. 
\begin{proof}[Proof of Proposition \ref{prop_diag}]
Recall that we have a commutative diagram:
\begin{equation}\label{diag1}
\xymatrix{
A\times_X B=A\cap B\ar[d]^{j_A}\ar[rr]^{\qquad j_B}&& A\ar[d]^{i_A}\\
B\ar[rr]^{i_B}&&X.
}
\end{equation}
And we need to prove the isomorphism (of derived functors)
\begin{equation}\label{eq_main}
i_{B!}j_{A*}(\sF|_{A\cap B})\cong i_{A*}j_{B!}(\sF|_{A\cap B}).
\end{equation}


To prove \eqref{eq_main}, we  consider another diagram, whose intersection with the diagonal gives rise to diagram \eqref{diag1}:
\begin{equation}\label{diag2}
\xymatrix{
A\times B\ar[d]^{i_A\times \id_B}\ar[rr]^{\id_A \times i_B}&& A\times X\ar[d]^{i_A\times \id_X}\\
X\times B\ar[rr]^{\id_X \times i_B}&&X\times X.
}
\end{equation}
We denote the diagonal map by $\Delta_X: X\to X\times X$, and denote its various restrictions by $\Delta_{A\cap B}: A\cap B\to A\times B$, $\Delta_{A}: A\to A\times X$, $\Delta_B: B\to X\times B$. These diagonal maps induce a morphism from the diagram \eqref{diag1} to the diagram \eqref{diag2}. 

With the above notations, we are ready to prove Proposition \ref{prop_diag}. By \cite[Proposition 10.2.9]{MS}, we have the following isomorphism, whose both sides are isomorphic to $\Delta_X^*(i_{A*}\sF \boxtimes  i_{B!}(\K_B))$: 
\begin{equation}\label{eq_iso}
\Delta_X^*(\id_X\times i_B)_!(i_{A}\times \id_B)_*(\sF \boxtimes  \K_B)\cong \Delta_X^*(i_{A}\times \id_X)_*(\id_A\times i_B)_!(\sF \boxtimes  \K_B).
\end{equation}
We can simplify the left-hand side of \eqref{eq_iso} as follows,
\begin{align*}
\Delta_X^*(\id_X\times i_B)_!(i_{A}\times \id_B)_*(\sF \boxtimes  \K_B)&\cong i_{B!}\Delta_B^*(i_{A}\times \id_B)_*(\sF \boxtimes  \K_B)\\
& \cong i_{B!}j_{A*}\Delta_{A\cap B}^*(\sF \boxtimes  \K_B)\\
&\cong i_{B!}j_{A*}(\sF|_{A\cap B})
\end{align*}
where the first isomorphism follows from 
the base change isomorphism for the direct image with proper support
(\cite[Theorem 2.3.26]{Dimca}), the second follows from \cite[Proposition 6.1.1]{Sch}, and the last is obvious. Let us check that  the transversality assumptions of \cite[Proposition 6.1.1]{Sch} hold here. Let $\mathfrak{S}$ be a Whitney stratification of $X$ such that $\sF$ is constructible with respect to $\mathfrak{S}$ and $A$ is a union of strata. Similarly, let $\mathfrak{T}$ be a Whitney stratification of $X$ such that $B_0$ is a union of strata. Now, $\mathfrak{S}\times g\cdot \mathfrak{T}$ is a Whitney stratification of $X\times X$ such that both $A\times B$ and $A\times X$ are unions of strata, and $\sF\boxtimes \K_B$ is constructible with respect to $\mathfrak{S}\times g\cdot \mathfrak{T}$. Thus, by Corollary \ref{cor_transversal}, the assumptions of \cite[Proposition 6.1.1]{Sch} hold. 
Also notice that in \cite[Proposition 6.1.1]{Sch}, the statement requires the corresponding embedding  $k'$ of loc, sit. to be an open embedding. Nevertheless, it also applies to the case when $k'$ is a locally closed embedding by factoring it as the open inclusion to the closure, which is also a union of strata, and the closed embedding to the ambient space, for which the needed isomorphism is just the 
base change isomorphism for the direct image with proper support.

Similarly, we can simplify the right-hand side of \eqref{eq_iso} by
\begin{align*}
\Delta_X^*(i_{A}\times \id_X)_*(\id_A\times i_B)_!(\sF \boxtimes  \K_B)&\cong i_{A*}\Delta_A^*(\id_A\times i_B)_!(\sF\boxtimes  \K_B)\\
&\cong i_{A*}j_{B!}\Delta_{A\cap B}^*(\sF\boxtimes  \K_B)\\
&\cong i_{A*}j_{B!}(\sF|_{A\cap B})
\end{align*}
where the first isomorphism follows from \cite[Proposition 6.1.1]{Sch}, the second follows from 
the base change isomorphism for the direct image with proper support,
and the last is obvious. By the above simplifications of both sides of \eqref{eq_iso}, we deduce the isomorphism~\eqref{eq_main}.
\end{proof}

\begin{proof}[Proof of Theorem \ref{thm_GV}]
Let $d=\codim_X B$, with $B$ smooth of pure codimension. Since $B$ is a generic translate and therefore transversal, $\sP|_{A\cap B}[-d]$ is a perverse sheaf (see e.g. \cite[Proposition 10.2.27]{MS} or \cite[Lemma 6.0.4]{Sch}). 
Then Proposition \ref{prop_diag} implies by the compactness of $X$:
\begin{align}\label{eq_align}
\begin{split}
H_c^i\big(B, j_{A*}(\sP|_{A\cap B}[-d])\big)&\cong H^i\big(X, i_{B!}j_{A*}(\sP|_{A\cap B}[-d])\big)\\
&\cong H^i\big(X, i_{A*}j_{B!}(\sP|_{A\cap B}[-d])\big)\\
&\cong H^i\big(A, j_{B!}(\sP|_{A\cap B}[-d])\big).
\end{split}
\end{align}
Since both $j_A$ and $j_B$ are quasi-finite and affine, both $j_{A*}(\sP|_{A\cap B}[-d])$ and $j_{B!}(\sP|_{A\cap B}[-d])$ are perverse sheaves 
(see e.g. \cite[Example 10.2.38, Theorem 10.3.69]{MS} and  \cite[Corollary 6.0.5, Theorem 6.0.4]{Sch}). Thus, by Artin's vanishing theorem 
(see e.g. \cite[Theorem 10.3.59]{MS} and \cite[Corollary 6.0.4]{Sch}), we have
\[
H_c^i\big(B, j_{A*}(\sP|_{A\cap B}[-d])\big)=0,\quad \text{for any}\enspace i< 0,
\]
and 
\[
H^i\big(A, j_{B!}(\sP|_{A\cap B}[-d])\big)=0,\quad  \text{for any}\enspace i> 0. 
\]
Combining with the isomorphism \eqref{eq_align}, we have 
\[
H_c^i\big(B, j_{A*}(\sP|_{A\cap B}[-d])\big) \cong H^i\big(A, j_{B!}(\sP|_{A\cap B}[-d])\big) =0,\quad \text{for any}\enspace i\neq 0,
\]
which is the desired vanishing after a shift of degree.
\end{proof}

\section{The Euler characteristics}
We prove Corollary \ref{cor_Euler} and Conjecture \ref{mainconj} in this section. 
\begin{proof}[Proof of Corollary \ref{cor_Euler}]
By Theorem \ref{thm_GV}, we have 
\[
(-1)^{\codim B}\chi_c\big(B, j_{A*}(\sP|_{A\cap B})\big)\geq 0.
\]
By \cite[Equation (10.3) and Example 10.4.37]{MS}  and \cite[Section 2.3
and Section 6.0.6]{Sch}), we have
\[
\chi_c\big(B, j_{A*}(\sP|_{A\cap B})\big)=\chi\big(B, j_{A*}(\sP|_{A\cap B})\big)=\chi(A\cap B, \sP|_{A\cap B}).
\]
Thus, the claimed inequality of Corollary \ref{cor_Euler} follows. 
\end{proof}

\begin{proof}[Proof of Conjecture \ref{mainconj}]
Let $A=X_{\lambda}^\circ\cap g X_\mu^\circ$ and $B_0=X_{\nu'}^\circ$. By Kleiman's theorem, $A$ is smooth pure dimensional, so that $\sP\coloneqq \bK_A[\dim A]$ is a perverse sheaf. By Corollary \ref{cor_Euler}, 
\[
(-1)^{\codim B}\chi(A\cap B, \sP|_{A\cap B})\geq 0.
\]
Since $\sP=\bK_A[\dim A]$, we have $\sP|_{A\cap B}=\bK_{A\cap B}[\dim A]$, and hence
\begin{align*}
(-1)^{\codim B}\chi(A\cap B, \sP|_{A\cap B})&=(-1)^{\dim A-\codim B}\chi(A\cap B)\\
&=(-1)^{\dim A\cap B}\chi(A\cap B)\\
&=(-1)^{d}\chi(X_\lambda^\circ\cap g X_\mu^\circ\cap h X_{\nu'}^\circ)\\
&=E_{\lambda, \mu, \nu'}.
\end{align*}
Therefore, we have $E_{\lambda, \mu, \nu'}\geq 0$. 
\end{proof}

Theorem \ref{thm_2} can be proved in the same way, as  explained before that Theorem.


\begin{thebibliography}{ADMSPRa}
\bibitem[AMSS17]{AMSS1} Paolo Aluffi, Leonardo C. Mihalcea, J\"{o}rg Sch\"{u}rmann, Changjian Su, \textit{Shadows of characteristic cycles, Verma modules, and positivity of Chern-Schwartz-MacPherson classes of Schubert cells}. arXiv:1708.08697.
To appear in Duke Math Journal.

\bibitem[AMSS22a]{AMSS0} Paolo Aluffi, Leonardo C. Mihalcea, J\"{o}rg Sch\"{u}rmann, Changjian Su, 
\textit{Positivity of Segre-MacPherson classes}. Facets of algebraic geometry. Vol. I, 1--28,
London Math. Soc. Lecture Note Ser., 472, Cambridge Univ. Press, Cambridge, 2022.

\bibitem[AMSS22b]{AMSS} Paolo Aluffi, Leonardo C. Mihalcea, J\"{o}rg Sch\"{u}rmann, Changjian Su, \textit{From motivic Chern classes of Schubert cells to their Hirzebruch and CSM classes}. arXiv:2212.12509. 

\bibitem[BC12]{BC} Sara Billey, Izzet Coskun, \textit{Singularities of Generalized Richardson Varieties}, Communications in Algebra 40
(2012), 1466--1495.

\bibitem[Bri12]{Bri12} Michel Brion, \textit{Spherical varieties}. In Highlights in Lie algebraic methods, volume 295 of Progr.
Math., pages 3--24. Birkh\"{a}user/Springer, New York, 2012.

\bibitem[Dim04]{Dimca} Alexandru Dimca,  \textit{Sheaves in topology}. Universitext. Springer-Verlag, Berlin, 2004.


\bibitem[FGX22]{FGX22} Neil Fan, Peter Guo, and Rui Xiong,  \textit{ Pieri and Murnaghan--Nakayama type rules for Chern
classes of Schubert cells}. arXiv:2211.06802.

\bibitem[FK00]{FK} Joseph Franecki, Mikhail Kapranov, 
\textit{The Gauss map and a noncompact Riemann-Roch formula for constructible sheaves on semiabelian varieties.}
Duke Math. J. 104 (2000), no. 1, 171--180.

\bibitem[GL96]{GaL} Ofer Gabber, François Loeser, 
\textit{Faisceaux pervers l-adiques sur un tore. }
Duke Math. J. 83 (1996), no. 3, 501--606.

\bibitem[GL87]{GL} Mark Green, Robert Lazarsfeld, \textit{Deformation theory, generic vanishing theorems, and some conjectures of Enriques, Catanese and Beauville.}
Invent. Math. 90 (1987), no. 2, 389--407.

\bibitem[Huh13]{Huh}  June Huh, \textit{The maximum likelihood degree of a very affine variety}. Compos. Math. 149  (2013), 1245--1266.

\bibitem[Kle74]{Kleiman} Steven Kleiman, \textit{The transversality of a general translate}. Compositio Math. 28 (1974), 287--297.

\bibitem[KZJ21]{KZJ} Allen Knutson, Paul Zinn-Justin, \textit{Schubert puzzles and integrability II: multiplying motivic Segre classes}. arXiv:2102.00563.

\bibitem[Knu22]{Kn22} Allen Knutson, \textit{Schubert calculus and quiver varieties}, available at\\ \url{https://pi.math.cornell.edu/~allenk/Knutson2022.pdf}

\bibitem[Kum22]{Kumar} Shrawan Kumar, \textit{Conjectural positivity of Chern-Schwartz-MacPherson classes for Richardson cells}. arXiv:2208.03527,
International Mathematics Research Notices, rnad036 (2023).


\bibitem[LMW19]{LMW} Yongqiang Liu, Laurentiu Maxim, Botong Wang, \textit{Generic vanishing for semi-abelian varieties and integral Alexander modules.} Math. Z. 293 (2019), no. 1-2, 629--645.

\bibitem[LMW21]{LMW2} Yongqiang Liu, Laurentiu Maxim, Botong Wang,  \textit{Euclidean distance degree of projective varieties}. 
Int. Math. Res. Not. IMRN Vo. 2021, no. 20, 15788--15802.

\bibitem[MS22]{MS} Laurentiu Maxim, J\"{o}rg Sch\"{u}rmann, \textit{Constructible sheaf complexes in complex geometry and applications}. Handbook of geometry and topology of singularities III, 679--791, Springer, Cham, 2022.


\bibitem[Ric92]{Ric92} R. W. Richardson,  \textit{Intersections of double cosets in algebraic groups}. Indag.
Math. (N.S.) 3 (1992), 69--77.

\bibitem[RW17]{RW} Jose Israel Rodriguez, Botong Wang, \textit{The maximum likelihood degree of mixtures of independence models}. SIAM J. Appl. Algebra Geom. 1 (2017), no. 1, 484--506.

\bibitem[Sch15]{Schnell} Christian Schnell, \textit{Holonomic D-modules on abelian varieties.} Publ. Math. Inst. Hautes Études Sci. 121 (2015), 1--55.

\bibitem[Sch03]{Sch} J\"{o}rg Sch\"{u}rmann, \textit{Topology of singular spaces and constructible sheaves}. Instytut Matematyczny Polskiej Akademii Nauk. Monografie Matematyczne (New Series), 63. Birkh\"{a}user Verlag, Basel, 2003.

\bibitem[Sch17]{Sch17}  J\"{o}rg Sch\"{u}rmann, \textit{Chern classes and transversality for singular spaces}. In Singularities in Geometry,
Topology, Foliations and Dynamics, Trends in Mathematics, pages 207--231. Birkh\"{a}user,
Basel, 2017.

\bibitem[ST10]{ST}  J\"{o}rg Sch\"{u}rmann and Mihai Tib\u{a}r. \textit{Index formula for MacPherson cycles of affine algebraic varieties}.
Tohoku Math. J.  62 (2010), 29--44.

\bibitem[STV05]{STV}  Jos\'{e} Seade, Mihai Tib\u{a}r, and Alberto Verjovsky. \textit{Global Euler obstruction and polar invariants}. Math. Ann., 333 (205), 393--403. 


\bibitem[Su21]{Su21} Changjian Su, \textit{Structure constants for Chern classes of Schubert cells}. Math. Z., 298 (2021), 193--213.

\end{thebibliography}
\end{document}